\documentclass[11pt]{amsart}
 \newtheorem{thm}{Theorem}[section]
 \newtheorem{cor}[thm]{Corollary}
 \newtheorem{lem}[thm]{Lemma}
 \newtheorem{prop}[thm]{Proposition}
 \theoremstyle{definition}
 \newtheorem{defn}[thm]{Definition}
 \theoremstyle{remark}
 
 \numberwithin{equation}{section}

 \newcommand{\norm}[1]{\left\Vert#1\right\Vert}
 
 \newcommand{\C}{\mathbb{C}}
 
\begin{document}

\title[]
 { Operator-Valued Norms}

\author{ Yun-Su Kim }


\email{kimys@indiana.edu}



\commby{Daniel J. Rudolph}


\begin{abstract}We introduce two kinds of operator-valued norms.
One of them is an $L(H)$-valued norm. The other one is an
$L(C(K))$-valued norm. We characterize the completeness with
respect to a bounded $L(H)$-valued norm. Furthermore, for a given
Banach space $\textbf{B}$, we provide an $L(C(K))$-valued norm on
$\textbf{B}$. and we introduce an $L(C(K))$-valued norm on a
Banach space satisfying special properties.

\end{abstract}

\maketitle

\section*{Introduction}
A nonnegative real-valued function $\norm{\emph{ }}$ defined on a
vector space, called  a \emph{norm}, is a very fundamental and
important project in analysis.\vskip0.1cm

Let $(X,\parallel$ $\parallel_{X})$ be a normed linear space which
is closed relative to the topology induced by the metric defined
by its norm. In this paper, $\textbf{B}$ and $H$ always denote a
Banach space and a Hilbert space, respectively. A lot of results
for scalar-valued functions have been extended to vector-valued
ones. As an example, in \cite{P}, W. L. Paschke introduced an
operator-valued inner product.

We introduce two kinds of operator-valued norms. In section
\ref{b2}, to introduce an \emph{$L(H)$-valued norm} $F$, a
function from $X$ to $L(H)$, we use positive operators instead of
nonnegative real numbers.

For any Banach space $\textbf{B}$, we have a trivial example
$F_{1}:\textbf{B}\rightarrow{L(H)}$ which is an $L(H)$-valued norm
defined by $F_{1}(b)=\norm{b}I_{H}$ where $I_{H}$ is the identity
mapping on $H$. Furthermore, with an $L(H)$-valued norm $F$,
automatically we introduce a notion of \emph{completeness with
respect to} $F$, and in Theorem \ref{b1} we  characterize the
completeness with respect to $F$.

Along the way we consider how to define an $L(\textbf{B})$-valued
norm on a normed linear space $X$. Since we do not have a notion
of positive operators in $L(\textbf{B})$, we have some
difficulties in defining $L(\textbf{B})$-valued norms. To
accomplish this, we use a fundamental Banach space $C(K)$ for some
compact Hausdorff space $K$.

 Thus, in section
\ref{a3}, we introduce a notion of an $L(C(K))$-valued norm,
instead of an $L(\textbf{B})$-valued norm, and obtain some
fundamental properties of $L(C(K))$-valued norms. By using
Banach's theorem, we provide an $L(C(K))$-valued norm on a given
Banach space $\textbf{B}$ in Theorem \ref{a6}.

\section{\textbf{$L(H)$-Valued Norms}}

\subsection{\textbf{L(H)-Valued Norms }}\label{b2}
Let $\C$ denote the set of complex numbers and $L(\textbf{B})$
denote the set of bounded operators on a Banach space
$\textbf{B}$. For $T\in{L(\textbf{B})}$, the norm of $T$ is
defined by $\norm{T}=\sup_{\norm{x}=1}\norm{Tx}$. Let
$(X,\parallel$ $\parallel_{X})$ be a normed linear space which is
closed relative to the topology induced by the metric defined by
its norm.

A nonnegative real-valued function $\norm{\emph{ }}$ defined on a
vector space, called  a \emph{norm}, is a very fundamental and
important project in analysis. To introduce the notion of
$L(H)$-valued norm on a normed linear space $X$ for a Hilbert
space $H$, we use positive operators.

Throughout this paper $H$ will denote a Hilbert space and for any
vectors $h_{1}$ and $h_{2}$ in $H$, $(h_{1},h_{2})$ is the inner
product of $h_{1}$ and $h_{2}$.

\begin{defn}\label{1}

(i) If a function $F:X\rightarrow{L(H)}$ has the following
properties :
\begin{enumerate}

\item For any $x\in{X}$, $F(x)\geq{0}$, i.e. $F(x)$ is a positive
operator.

\item (Triangle Inequality) $F(x+y)\leq{F(x)+F(y)}$ for any $x$
and $y$ in $X$.

\item $F(\lambda{x})=|\lambda|F(x)$ for any $\lambda\in{\C}$.

\item $F(x)=0$ if and only if $x=0$, \end{enumerate} then $F$ is
said to be an $L(H)$-\emph{valued norm } defined on $X$.

\noindent (ii) If $\sup_{\norm{x}=1}\norm{F(x)}<\infty$, then $F$
is said to be \emph{bounded}.

\end{defn}

Let $\mathbb{T}$ be the unit circle in the complex plane. As an
example of such projects, define
$F:L^{\infty}(\mathbb{T})\rightarrow{L(L^{2}(\mathbb{T}))}$ by
\[F(g)=M_{g},\]
where $L^{p}(\mathbb{T})$ is the Lebesgue space with respect to
the Lebesgue measure $\mu$ on $\mathbb{T}$ such that
$\mu(\mathbb{T})=1$, and
$M_{g}:L^{2}(\mathbb{T})\rightarrow{L^{2}(\mathbb{T})}$ is the
bounded operator defined by $M_{g}(f)=|g|\cdot{f}$ for
$f\in{L^{2}(\mathbb{T})}$. Then clearly, $F$ is an
$L(L^{2}(\mathbb{T}))$-valued norm on $L^{\infty}(\mathbb{T})$.




\begin{lem}
Let $F:H\rightarrow{L(H)}$ be an $L(H)$-valued norm on $H$ and
$T\in{L(H)}$ be an injective operator.

Then $F\circ{T}$ is an $L(H)$-valued norm on $H$.

\end{lem}
\begin{proof}
Let $G=F\circ{T}$ and $h\in{H}$ be given. Since $F(Th)$ is a
positive operator, $F(Th)=S^{\ast}S$ for some $S\in{L(H)}$. Then
\begin{center}$(G(h)k,k)=\norm{Sk}^{2}\geq{0}$ for $k\in{K}$.\end{center} Thus $G(h)\geq{0}$ for any
$h\in{H}$.

Since $T$ is linear and $F$ is an $L(H)$-valued norm on $H$,
$G(h+k)=F(Th+Tk)\leq{F(Th)+F(Tk)}={G(h)+G(k)}$ for $h$ and $k$ in
$H$, and $G(\lambda{h})=F(\lambda{Th})=|\lambda|G(h)$ for any
$\lambda\in{\C}$.

Let $G(x)=0$ for some $x\in{H}$. Then $F(Tx)=0$ and so $Tx=0$,
since $F$ is an $L(H)$-valued norm. The injectivity of $T$ implies
that $x=0$. Conversely, if $x=0$, then clearly $G(x)=0$.
Therefore, $G=F\circ{T}$ is an $L(H)$-valued norm.

\end{proof}

\begin{prop}\label{6}\cite{R}
If $T\in{L(H)}$ is normal, then
\begin{center}{$\norm{T}=\sup\{|(Th,h)|:h\in{H},
\norm{h}\leq{1}\}.$}\end{center}
\end{prop}

\begin{prop}\label{5}
Let $F:X\rightarrow{L(H)}$ be an $L(H)$-valued norm on $X$. Then
\begin{enumerate}
\item $\norm{F(x+y)}\leq\norm{F(x)}+\norm{F(y)}$ for any $x$ and
$y$ in $X$.

\item $\norm{F(x)-F(y)}\leq{\norm{F(x-y)}}$ for any $x$ and $y$ in
$X$.\end{enumerate}
\end{prop}
\begin{proof}
(1) For any $x$ and $y$ in $X$, by triangle inequality,
$0\leq{F(x+y)}\leq{F(x)+F(y)} $. It follows that
\begin{center}$\norm{F(x+y)}\leq\norm{F(x)+F(y)}\leq\norm{F(x)}+\norm{F(y)}.$\end{center}
\vskip0.2cm

(2) By triangle inequality, for any $x$ and $y$ in $X$,
\begin{center}$F(x)\leq{F(x-y)+F(y)}$ or $F(x)-F(y)\leq{F(x-y)}$.\end{center} Thus for
any $h\in{H}$,
\begin{equation}\label{2}
((F(x)-F(y))h,h)\leq{(F(x-y)h,h)}. \end{equation} Similarly, we
have $F(y)-F(x)\leq{F(y-x)}=|-1|F(x-y)=F(x-y)$. Thus for any
$h\in{H}$,
\begin{equation}\label{3}
((F(y)-F(x))h,h)\leq{(F(x-y)h,h)}.
\end{equation}
Inequalities (\ref{2}) and (\ref{3}) imply that
\begin{equation}\label{4}
|((F(x)-F(y))h,h)|\leq|(F(x-y)h,h)|.
\end{equation}

Since $F(x)$ and $F(y)$ are positive operators, from Proposition
\ref{6} and inequality (\ref{4}), we conclude that
\begin{center}$\norm{F(x)-F(y)}\leq\norm{F(x-y)}.$\end{center}
\end{proof}
\begin{prop}\label{7}
Let $F:X\rightarrow{L(H)}$ be an $L(H)$-valued norm on $X$. Then
$F$ is continuous at 0 if and only if $F$ is a continuous function
on $X$.
\end{prop}
\begin{proof}
Suppose that $F$ is continuous at 0. Let $x\in{X}$ and
$\{x_{n}\}_{n=1}^{\infty}$ be a sequence in $X$ such that
\begin{equation}\label{a1}\norm{x_{n}-x}_{X}\rightarrow{0}\end{equation} as
$n\rightarrow\infty$. By Proposition \ref{5},
\begin{equation}\label{a}
\norm{F(x_{n})-F(x)}\leq\norm{F(x_{n}-x)}.\end{equation}

Since $F$ is continuous at 0 and $F(0)=0$, (\ref{a1}) and
(\ref{a}) imply that
\begin{center}$\norm{F(x_{n})-F(x)}\rightarrow{0}$\end{center}
as $n\rightarrow\infty$. Thus $F$ is continuous at $x\in{X}$.
\vskip0.2cm The converse is clear.
\end{proof}
\subsection{\textbf{Completeness with respect to $F$}}\label{b}
In Section \ref{b2}, we provided new operator-valued norm and so
in this section we provide a notion of completeness with respect
to an $L(H)$-valued norm.
\begin{defn}\label{b}
Let $\{x_{n}\}_{n=1}^{\infty}$ be a Cauchy sequence in $X$ and
$F:X\rightarrow{L(H)}$ be an $L(H)$-valued norm on $X$.

If $\{F(x_{n})\}_{n=1}^{\infty}$ is also a Cauchy sequence in
$L(H)$, then $X$ is said to be \emph{complete with respect to} $F$
or \emph{a Banach space with respect to }$F$.
\end{defn}

\begin{thm}\label{b1}
Let $F:X\rightarrow{L(H)}$ be a bounded $L(H)$-valued norm on $X$.
$X$ is complete with respect to $F$ if and only if $F$ is
continuous at 0.
\end{thm}
\begin{proof}
Suppose $X$ is complete with respect to $F$ and
$\{x_{n}\}_{n=1}^{\infty}$ is a sequence in $X$ such that
\begin{equation}\label{s1}\lim_{n \rightarrow \infty}{x_{n}}=0.\end{equation} Then $\{x_{n}\}_{n=1}^{\infty}$ is a Cauchy
sequence in $X$. Since $X$ is complete with respect to $F$,
$\{F(x_{n})\}_{n=1}^{\infty}$ is also a Cauchy sequence.
Since $L(H)$ is a Banach
space, there is an operator $T\in{L(H)}$ such that
\begin{equation}\label{s}\lim_{n \rightarrow \infty}{F(x_{n})}=T.\end{equation}

If $A=\{n:x_{n}=0\}$ is infinite, then for $n_{k}\in{A}$
$(k=1,2,3,\cdot\cdot\cdot)$, \begin{center}$\lim_{k \rightarrow
\infty}{F(x_{n_{k}})}=F(0)=0$\end{center} and so by equation
(\ref{s}), $T=0$. \vskip0.2cm If $A=\{n:x_{n}=0\}$ is not
infinite, then there is a nonzero subsequence
$\{y_{n}\}_{n=1}^{\infty}$ of $\{x_{n}\}_{n=1}^{\infty}$ such that
\begin{equation}\label{s2}\lim_{n \rightarrow \infty}{y_{n}}=0\texttt{ and }\lim_{n \rightarrow \infty}{F(y_{n})}=T.
\end{equation}
By equation (\ref{s2}),
\begin{equation}\label{s3}\lim_{n \rightarrow
\infty}\norm{F(y_{n})}=\norm{T}.\end{equation} Let
$M_{F}=\sup_{\norm{x}=1}{\norm{F(x)}}(<\infty)$. By property (3)
of Definition \ref{1} and \ref{s3}, $\norm{T}=\lim_{n \rightarrow
\infty}\norm{y_{n}}_{X}\norm{F(\frac{y_{n}}{\norm{y_{n}}_{X}})}\leq\lim_{n\rightarrow
\infty}\norm{y_{n}}_{X}\cdot{M_F}$.

It follows that \begin{center}$\norm{T}=0$\end{center} and so
$T=0$. Thus $T\equiv{0}$ whether $A$ is infinite or not.

\noindent By (\ref{s}), $\lim_{n \rightarrow
\infty}{F(x_{n})}=0=F(0)$ which proves that $F$ is continuous at
0. \vskip0.2cm

Conversely, suppose that $F$ is continuous at 0. By Proposition
\ref{5} (2), it is easy to see that $X$ is complete with respect
to $F$.
\end{proof}
From Proposition \ref{7}, we obtain the following result:
\begin{cor}
Let $F:X\rightarrow{L(H)}$ be a bounded $L(H)$-valued norm on $X$.
$X$ is complete with respect to $F$ if and only if $F$ is
continuous on $X$.
\end{cor}

\section{\textbf{$L(C(K))$-Valued Norms}}

Let $K$ be a compact Hausdorff space. In this section, we
introduce a new norm whose value is a bounded operator on a Banach
space and provide some examples of that.
\subsection{\textbf{$L(C(K))$-Valued Norms}}\label{a3}
Many of well-known function spaces are Banach spaces. Thus we
consider how to define $L(\textbf{B})$-valued norms. Since we do
not have the notion of positive operators in $L(\textbf{B})$, we
can not provide the same definition as \ref{1}. To accomplish our
results, we provide an $L(C(K))$-valued norm, instead of an
$L(\textbf{B})$-valued norm.



\begin{defn}\label{1}
If a function $F:X\rightarrow{L(C(K))}$ has the following
properties :
\begin{enumerate}

\item For any $x\in{X}$, $(F(x))(f)\geq{0}$ whenever $f\geq{0}$.

\item (Triangle Inequality) If $f\geq{0}$, then
${(F(x)+F(y)-F(x+y))f}\geq{0}$ for any $x$ and $y$ in $X$.

\item $F(\lambda{x})=|\lambda|F(x)$ for any $\lambda\in{\C}$.

\item $F(x)=0$ if and only if $x=0$, \end{enumerate} then $F$ is
said to be an $L(C(K))$\emph{-valued norm} defined on $X$.
\end{defn}

Since $F(x)\in{L(C(K))}$ for any $x\in{X}$,
\begin{center}\[\norm{F(x)}=\sup_{f(\neq{0})\in{C(K)}}\frac{\norm{F(x)f}}{\norm{f}}<\infty.\]\end{center}
If \begin{center}$\sup_{\norm{x}=1}{\norm{F(x)}}$
$<\infty$,\end{center} then $F$ is said to be \emph{bounded}.
\vskip0.2cm As an example, let $X=H^{\infty}[0,1]$ which is the
space of all analytic bounded Lebesgue measurable functions $g$ on
$[0,1]$ with a norm defined by
\begin{center}$\norm{g}_{\infty}=\texttt{essential supremum of }|g|,$\end{center}
and $C[0,1]$ be the space of all continuous functions $f$ on
$[0,1]$. For each $g\in{H^{\infty}[0,1]}$, we define a
multiplication operator $M_{g}$ on $C[0,1]$ by
\begin{equation}\label{34}M_{g}(f)=|g|\cdot{f}.\end{equation} Define
$F:H^{\infty}[0,1]\rightarrow{L(C[0,1])}$ by
\begin{center}$F(g)=M_{g}$.\end{center} We can easily see that $F$ is an
operator-valued norm on $H^{\infty}[0,1]$ and
\begin{center}$\sup_{\norm{g}_{\infty}=1}{\norm{F(g)}}=
\sup_{\norm{g}_{\infty}=1}{\norm{M_{g}}}=\sup_{\norm{g}_{\infty}=1}\norm{g}_{\infty}=1$,\end{center}
that is, $F$ is bounded.

\begin{defn}\cite{D}
Let $\textbf{B}$ be a Banach algebra. A complex linear functional
$\varphi$ on $\textbf{B}$ is said to be \emph{multiplicative} if :
\begin{enumerate}
\item $\varphi(fg)=\varphi(f)\varphi(g)$ for $f$ and $g$ in
$\textbf{B}$; and \item $\varphi(1)=1.$\end{enumerate}

\end{defn}
\noindent The set of all multiplicative linear functionals on
$\textbf{B}$ is denoted by $\textbf{M}_{\textbf{B}}$.
\begin{defn}\cite{D}
For the Banach algebra $\textbf{B}$, if
$\textbf{M}_{\textbf{B}}\neq\emptyset$, then the Gelfand transform
is the function
$\Gamma:\textbf{B}\rightarrow{C(\textbf{M}_{\textbf{B}})}$ given
by
\begin{center}$\Gamma(f)(\varphi)=\varphi(f)$\end{center}
for $\varphi$ in $\textbf{M}_{\textbf{B}}$.
\end{defn}
\begin{prop}\cite{D}\label{a7}
If $\textbf{B}$ is a Banach algebra and $\Gamma$ is the Gelfand
transform on $\textbf{B}$, then
\begin{enumerate}
\item $\Gamma$ is an algebra homomorphism; and

\item $\norm{\Gamma{f}}_{\infty}\leq\norm{f}$ for $f$ in
$\textbf{B}$.
\end{enumerate}
\end{prop}
\begin{prop}\cite{D}\label{a8}
Let $Y$ be a compact Hausdorff space. If $\textbf{B}$ is a closed
self-adjoint subalgebra of $C(Y)$ containing the constant function
1, then the Gelfand transform $\Gamma$ is an isometric isomorphism
from $\textbf{B}$ onto $C(\textbf{M}_{\textbf{B}})$.
\end{prop}

\begin{cor}\label{a9}
Let $Y$ be a compact Hausdorff space. If $\textbf{B}$ is a closed
self-adjoint subalgebra of $C(Y)$ containing the constant function
1, then there is a bounded $L(C(\textbf{M}_{\textbf{B}}))$-valued
norm
\begin{center}$F:\textbf{B}\rightarrow{L(C(\textbf{M}_{\textbf{B}}))}$\end{center} such that
\begin{center}$F(ab)=F(a)F(b)$\end{center} for any $a$ and $b$ in $\textbf{B}$.
\end{cor}
\begin{proof}
Define $F:\textbf{B}\rightarrow{L(C(\textbf{M}_{\textbf{B}}))}$ in
a similar way as Theorem \ref{a6} by using the Gelfand transform,
for $b$ in $\textbf{B}$,

\begin{center}$F(b)=M_{|\Gamma{b}|}$\end{center}
where
$M_{|\Gamma{b}|}:C(\textbf{M}_{\textbf{B}})\rightarrow{C(\textbf{M}_{\textbf{B}})}$
is defined by $M_{|\Gamma{b}|}(\varphi)=|\Gamma{b}|\cdot{\varphi}$
which means that
$(M_{|\Gamma{b}|}\varphi)(f)=|(\Gamma{b})f|\cdot{\varphi(f)}$ for
$\varphi\in{C(\textbf{M}_{\textbf{B}})}$ and
$f\in{\textbf{M}_{\textbf{B}}}$. Then $F$ is well-defined.

Let $a$ and $b$ be in $\textbf{B}$. Clearly, $F(b)\varphi\geq{0}$
whenever $\varphi\geq{0}$. If $\varphi\geq{0}$ and
$f\in{\textbf{M}_{\textbf{B}}}$, then
\[[(F(a)+F(b)-F(a+b))\varphi](f)=(|{(\Gamma{a})f}|+|{(\Gamma{b})f}|-|{(\Gamma{a+b})f}|)\varphi(f)\]

\begin{center}  $\quad\quad\quad\quad\quad\quad\quad\quad\quad\quad\quad\quad\quad\quad=(|f(a)|+|f(b)|-|f(a+b)|)\varphi(f)\geq{0}$\end{center}
which implies the triangle inequality.

Next, $F(\lambda{b})=|\lambda|F(b)$ for $\lambda$ in $\C$.

Finally, if $F(a)=0$, then $|\Gamma(a)|\cdot\varphi=0$ for any
$\varphi$ in $C(\textbf{M}_{\textbf{B}})$. If $\varphi\equiv{1}$,
then we have $|\Gamma(a)|=0$. Since $\Gamma$ is isometric by
Proposition \ref{a8}, $a=0$. Conversely, if $a=0$, then $F(a)=0$
is clear. Therefore, $F$ is an
$L(C(\textbf{M}_{\textbf{B}}))$-valued norm on $B$.

For any $b\in{\textbf{B}}$,
\begin{center}$\norm{F(b)}=\norm{M_{|\Gamma{b}|}}=\norm{\Gamma{b}}_{\infty}=\norm{b}$.\end{center}
Thus $F$ is bounded.

By Proposition \ref{a7}, $\Gamma$ is an algebra homomorphism. It
follows that for any $a$ and $b$ in $B$,
\begin{center}$F(ab)=M_{|\Gamma{(ab)}|}=M_{|\Gamma(a)\Gamma(b)|}
=M_{|\Gamma(a)|}M_{|\Gamma(b)|}=F(a)F(b)$\end{center} which proves
this Corollary.

\end{proof}

\begin{prop}\textbf{(Banach)} \cite{D}\label{s7}
Every Banach space $\textbf{B}$ is isometrically isomorphic to a
closed subspace of $C(K)$ where $K=(\textbf{B}^{*})_{1}$ is the
unit ball of the dual of the Banach space $\textbf{B}$ and $K$ is
given the weak$^{\ast}$-topology.
\end{prop}

\begin{thm}\label{a6}
Let $K$ be defined in the same way as Proposition \ref{s7}. For
any Banach space $(\textbf{B},\norm{\cdot}_{\textbf{B}})$, there
is a function $F:\textbf{B}\rightarrow{L(C(K))}$ such that $F$ is
a bounded operator-valued norm on $\textbf{B}$. In particular,
\begin{center}$\norm{F(b)}=\norm{b}$\end{center} for any
$b\in{\textbf{B}}$.
\end{thm}
\begin{proof}By Proposition \ref{s7}, there is a compact Hausdorff
space $K$ such that $\beta:\textbf{B}\rightarrow{C(K)}$ is a
linear map and isometric.

Define a function $F:\textbf{B}\rightarrow{L(C(K))}$ by
\begin{center}
$F(b)(f)=M_{|\beta{b}|}(f)=|\beta{b}|\cdot{f}$\end{center} for
$b\in{\textbf{B}}$ and $f\in{C(K)}$. Clearly, $F$ is well defined
and if $f\geq{0}$, then
\begin{equation}\label{8}F(b)f\geq{0}\end{equation} for any $b\in{\textbf{B}}$, and
\begin{equation}\label{9}(F(a)+F(b)-F(a+b))f=({M_{|\beta{a}|+|\beta{b}|}}-M_{|\beta(a+b)|})f\geq{0}
\end{equation}
for any $a\in{\textbf{B}}$ and $b\in{\textbf{B}}$.

For any $\lambda\in{\C}$ and $b\in{\textbf{B}}$,
\begin{equation}\label{10}F(\lambda{b})=M_{|\beta(\lambda{b})|}=M_{|\lambda||\beta{b}|}=|\lambda|M_{|\beta{b}|}=
|\lambda|F(b).\end{equation}

Finally, $F(b)=0$ if and only if $M_{|\beta{b}|}=0$ if and only if
$|\beta{b}|=0$ if and only if $b=0$, since $\beta$ is isometric.
Thus \begin{equation}\label{11}F(b)=0\texttt{ if and only if }
b=0\end{equation} for any $b\in{\textbf{B}}.$ From (\ref{8}),
(\ref{9}), (\ref{10}), and (\ref{11}), we conclude that $F$ is an
$L(C(K))$-valued norm on $\textbf{B}$.

Since for any $f\in{C(K)}$, $\norm{M_{f}}=\norm{f}_{\infty}$, and
$\beta$ is isometric, for any $b\in{\textbf{B}}$,
\begin{center}$\norm{F(b)}=\norm{M_{|\beta{b}|}}=\norm{\beta{b}}_{\infty}=\norm{b}$\end{center}
which proves this theorem.
\end{proof}

------------------------------------------------------------------------

\end{document}